\documentclass[11pt]{article}
\usepackage[a4paper,left=25mm,right=25mm,bottom=25mm,top=25mm]{geometry}
\usepackage[T2A]{fontenc}
\usepackage[cp1251]{inputenc}
\usepackage[russian]{babel}
\usepackage{amsthm}
\usepackage{amssymb}
\usepackage{amsmath}
\usepackage{amsfonts}
\usepackage{indentfirst}
\usepackage{multirow}
\usepackage{graphicx}
\usepackage{wrapfig}
\usepackage{titlesec}
\titleformat{\section}[block]{\large\bfseries\filcenter}{\large\bfseries\thesection.\,}{0pt}{}
\titleformat{\subsection}[runin]{\normalfont\bfseries}{\indent\normalfont\bfseries\thesubsection.\,}{0pt}{}[.]
\newtheoremstyle{theorem}%  % name
{0pt}%                      % space above
{0pt}%                      % space below
{\it}%                      % body font
{\parindent}%               % indent amount
{\bf}%                      % theorem head font
{\normalfont{\bf .}}%           % punctuation after theorem head
{ }%                        % space after theorem head
{}%                         % theorem head spec
\theoremstyle{theorem}
\newtheorem{theorem}{Theorem}[section]
\newtheorem{corollary}[theorem]{Corollary}

\newtheorem{lemma}[theorem]{Lemma}

\newtheoremstyle{definition}%  % name
{0pt}%                      % space above
{0pt}%                      % space below
{}%                         % body font
{\parindent}%               % indent amount
{\bf}%                      % theorem head font
{\normalfont{\bf .}}%           % punctuation after theorem head
{ }%                        % space after theorem head
{}%                         % theorem head spec
\theoremstyle{definition}
\theoremstyle{definition}
\newtheorem{remark}[theorem]{Remark}
\newtheorem{example}[theorem]{Example}

\makeatletter
\renewcommand\@biblabel[1]{#1.}
\makeatother
\begin{document}

\parindent=10mm
\begin{center}
\large\bfseries POLYNOMIALS OVER DIVISION RINGS
\end{center}

\medskip

\begin{center}
\large\bfseries A.~G.~Goutor$^1$, S.~V.~Tikhonov$^2$
\end{center}

\medskip

\begin{center}

$^1$ \it Belarusian State Iniversity, Minsk, Belarus \\
$^2$ \it Belarusian State Iniversity, Minsk, Belarus\\
e-mail: goutor7@gmail.com, tikhonovsv@bsu.by,
\end{center}

We consider properties of polynomials with coefficients in division rings.

A theorem on the decomposition of a polynomial with coefficients in an arbitrary division ring is obtained.
It is shown that if a non-central element is not a root of a polynomial over an arbitrary division ring,
then the conjugacy class of this element contains infinitely many elements that are not roots of this polynomial.
The paper also contains estimates for the number of different conjugacy classes of spherical roots for some types of polynomials over quaternion division algebras.


%%%%%%%%%%%%%%%%%%%%%%%%%%%%%%%%%%%%%%%%%%%%%%%%%%%%%%%%%%%%%%%%%%%%%%%%%%%%%%%%%%%%%%%%%%%%%%%%%%%%%%%%%%%%%%%%%%%%%%%%%%%%%%%%%%%%%%%%%%%

\section{Introduction and preliminary results}

Let $\cal R$ be a non-commutative associative division ring, $\cal R^*$ its group of units. ${\cal R}[x]$ denotes the ring of polynomials in $x$ with coefficients in $\cal R$,
we assume that $x$ commutes with elements of $\cal R$. Thus, every polynomial in ${\cal R}[x]$ has the form

\begin{equation} \label{eq_P(x)}
P(x)=a_nx^n+a_{n-1}x^{n-1}+\dots + a_1x+a_0,\quad a_0, \dots, a_n \in {\cal R}.
\end{equation}

Addition and multiplication of polynomials in ${\cal R}[x]$ are defined in a natural way.
The degree of a polynomial of the form (\ref{eq_P(x)}) is also defined in the usual way and is equal to $n$ if $a_n\ne 0$.
In the ring ${\cal R}[x]$ a right-hand division with remainder and a right-hand Euclidean algorithm hold, and for polynomials
$P(x),S(x) \in {\cal R}[x]$ their greatest common right divisor $\mbox{GCRD}(P(x),S(x))$ is defined (see \cite{Or33}).

The main properties of polynomials over division rings are described in \cite[Ch. 5, \S 16]{La91} (see also \cite{GoMo65}, \cite{BrWh83}).

The following theorem is proved in \cite{Be79}.

\begin{theorem} \label{th:Beck}
Let $\cal Q$ be a division quaternion algebra over the field $K$. Then every polyno\-mi\-al $P(x) \in {\cal Q}[x]$ can be
uniquely expressed as a product $P(x)=cG(x)H(x)$, where $c \in {\cal Q}^*$ is the leading coefficient of $P(x)$,
$H(x)$ is a monic polynomial with coefficients in $K$, and $G(x) \in {\cal Q}[x]$ is a monic polynomial
not divisible on the right by any non-constant polynomial from $K[x]$. Moreover, if ${\cal Q}[x]$ is treated as a free module of rank 4 over $K[x]$
with the standard basis $1,i,j,k$, then $H(x)$ is the greatest common divisor (in $K[x]$) of the coordinates of $P(x)$ in this basis.
\end{theorem}

One of the goals of this paper is to generalize this theorem to any division ring (see Theorem \ref{th:factorization} below).

For $a\in R$, define $P(a)$ as an element

$$
P(a)=a_na^n+a_{n-1}a^{n-1}+\dots + a_1a + a_0.
$$

We call an element $a\in R$ a (right) root of $P(x)$ if $P(a)=0$. It is known that $a\in R$ is a root of the polynomial $P(x)$ if and only if
$x-a$ is a right divisor of $P(x)$ in ${\cal R}[x]$ (\cite[Proposition 16.2]{La91}, i.e. $P(x)=H(x)(x-a)$ for some polynomial $H(x)$ in ${\cal R}[x]$.
Note that the equality $P(x)=H(x)S(x) \in {\cal R}[x]$ does not imply that $P(a)=H(a)S(a)$. In particular, if $a$ is a root of the polynomial $H(x)$,
then $a$ may not be a root of the polynomial $P(x)$.

The conjugacy class of an element $a\in {\cal R}$, which we will denote by $[a]$, consists of all elements of the form $qaq^{-1}$, where $q$ is an arbitrary nonzero element of $\cal R$.
The ring $\cal R$ splits into disjoint conjugacy classes. By ${\cal R}^{(c)}$ we denote the set of all elements of $\cal R$ that commute with an element $c \in {\cal R}$. ${\cal R}^{(c)}$
is a division subring of $\cal R$.

In the case of polynomials over fields, every polynomial of degree $n$ has at most $n$ roots. In the case of polynomials with coefficients in division rings, the situation is different:
a polynomial of degree $n$ can have infinitely many roots. The Gordon-Motzkin theorem (\cite[Th. 16.4]{La91}) says that a polynomial of degree $n$ in ${\cal R}[x]$ can have
roots in at most $n$ conjugacy classes of the ring $\cal R$. Moreover, if $P(x)\in {\cal R}[x]$ has two distinct roots in a conjugacy class, then $P(x)$ has
infinitely many roots in this class (see \cite[Theorem 16.11]{La91} and \cite[Proposition 3]{BrWh83}). In the case of a division quaternion algebra $\cal Q$,
if $P(x)\in {\cal Q}[x]$ has two distinct roots in a conjugacy class, then every element of this class is a root of $P(x)$.
This means that polynomials over division quaternion algebras have only two types of roots:
isolated and spherical. A root $q$ of a polynomial $P(x)$ is called spherical if $q$ does not belong to the center of the algebra and any element $d \in [q]$ is also a root of $P(x)$.
A root $q$ is called isolated if the conjugacy class $[q]$ contains only one root of the polynomial $P(x)$.
However, in the case where the minimal polynomial of the conjugacy class has degree greater than two, the situation is fundamentally different. In \cite{GoTi24}, for any conjugacy class with
the minimal polynomial of degree > 2, a quadratic polynomial is constructed that has infinitely many roots in this  conjugacy class, and there are infinitely many
elements in this conjugacy class that are not roots of this polynomial. In this paper, we show that if a non-central element $c$ is not a root of a polynomial over an arbitrary division ring,
then there are infinitely many elements in the
class $[c]$ that are not roots of this polynomial (see Theorem \ref{th:roots} below).
Also in this paper, estimates are obtained for the number of different conjugacy classes of spherical roots for some polynomials over quaternion algebras.

%%%%%%%%%%%%%%%%%%%%%%%%%%%%%%%%%%%%%%%%%%%%%%%%%%%%%%%%%%%%%%%%%%%%%%%%%%%%%%%%%%%%%%%%%%%%%%%%%%%%%%%%%%%%%%%%%%%%%%%%%%%%%%%%%%%%%%%%%%%%%%%%%%%%%%
%%%%%%%%%%%%%%%%%%%%%%%%%%%%%%%%%%%%%%%%%%%%%%%%%%%%%%%%%%%%%%%%%%%%%%%%%%%%%%%%%%%%%%%%%%%%%%%%%%%%%%%%%%%%%%%%%%%%%%%%%%%%%%%%%%%%%%%%%%%%%%%%%%%%%%

\section{Polynomials over arbitrary division rings}

Let ${\cal R}_0$ be a division subring of the ring ${\cal R}$, $\{c_i\}_{i\in I}$ a basis of the right vector space ${\cal R}$ over ${\cal R}_0$.
Then every polynomial $P(x) \in {\cal R}[x]$ can be uniquely represented as
\begin{equation} \label{eq_R0}
P(x) = \sum_{i \in I} c_i b_i(x),
\end{equation}
where $b_i(x) \in {\cal R}_0[x]$ are almost all zero.
The polynomials $b_i(x)$ are obtained as follows. Consider a polynomial $P(x)$ of the form (\ref{eq_P(x)}). Let us expand each coefficient $a_i, i=1,..., n$, of the polynomial in the basis $\{c_i\}_{i\in I}$. Then
$$
P(x)=(\sum_{i \in I} c_i a_{n,i})x^n+(\sum_{i \in I} c_i a_{n-1,i})x^{n-1}+\dots
+(\sum_{i \in I} c_i a_{1,i})x+(\sum_{i \in I} c_i a_{0,i})=
$$
$$
\sum_{i \in I} c_i (a_{n,i}x^n+a_{n-1,i}x^{n-1}+...+a_{1,i}x+a_{0,i}),
$$
where $a_{k,i} \in {\cal R}_0, k=0,..., n, i\in I$.

Thus, $b_i(x)=a_{n,i}x^n+a_{n-1,i}x^{n-1}+...+a_{1,i}x+a_{0,i}$ in the expansion (\ref{eq_R0}). Note that almost all of these polynomials are 0.
In the notation above, the following theorem holds.

\begin{theorem} \label{th:R0}
1. A polynomial $h(x)\in {\cal R}_0[x]$ is a right divisor of the polynomial $P(x)\in {\cal R}[x]$ if and only if $h(x)$ is a common right divisor of the polynomials $b_i(x)$, $i\in I$.

2. A polynomial $h(x)\in {\cal R}_0[x]$ is a right divisor of the greatest degree in ${\cal R}_0[x]$ of the polynomial $P(x)\in {\cal R}[x]$ if and only if
$h(x) = \mbox{GCRD}(b_i(x), i\in I)$ in ${\cal R}_0[x]$.

3. An element $\alpha \in {\cal R}_0$ is a root of the polynomial $P(x)\in {\cal R}[x]$ if and only if $\alpha$ is a common root of the polynomials $b_i(x), i\in I$.

4. If a polynomial $h(x)\in {\cal R}_0[x]$ is a right divisor of the greatest degree in ${\cal R}_0[x]$ of the polynomial $P(x)\in {\cal R}[x]$, then any root in ${\cal R}_0$
of the polynomial $P(x)$ is a root of the polynomial $h(x)$.

\end{theorem}

Proof. 1. If $h(x)$ divides every polynomial $b_i(x)$ on the right, then $h(x)$ divides $P(x)$ on the right as well.

Now let $P(x)=G(x)h(x)$ for some polynomial $G(x)\in {\cal R}[x]$. For the polynomial $G(x)$ there is a decomposition of the form (\ref{eq_R0}):
$$
G(x)=\sum_{i \in I} c_i d_i(x) ,
$$
where $d_i(x)$ is from ${\cal R}_0[x]$.
Then
$$
P(x)=\left(\sum_{i \in I} c_i d_i(x)\right)h(x)=\sum_{i \in I} c_i d_i(x)h(x).
$$
Thus, $b_i(x)=d_i(x)h(x)$ for all $i\in I$. Therefore, $h(x)$ is a common divisor of the polynomials $b_i(x)$, $i \in I$.

2. Follows from 1.

3. Follows from 1, since $\alpha$ is a root of the polynomial $P(x)$ if and only if $x-\alpha$ is a right divisor of the polynomial $P(x)$.

4. Follows from 2 and 3. Indeed, if $\alpha \in R_0$ is a root of the polynomial $P(x)$, then $x-\alpha$ is a common right divisor of the polynomials $b_i(x)$, $i \in I$. Then $x-\alpha$ is a
right divisor of their greatest common right divisor.  \qed

%%%%%%%%%%%%%%%%%%%%%%%%%%%%%%%%%%%%%%%%%%%%%%%%%%%%%%%%%%%%%%%%%%%%%%%%%%%%%%%%%%%%%%%%%%%%%%%%%%%%%%%%%%

\begin{remark}
If we consider a division algebra under the conditions of the previous theorem, then all roots of the polynomial $P(x)$ can be sought as roots of polynomials from subfields of the algebra
(for example, a root $a$ lies in the subfield $F(a)$). Thus, the problem of finding roots of a polynomial with coefficients in some algebra is reduced to the problem of finding roots in subfields.
\end{remark}

\begin{remark}
Let us show that the converse of item 4 of Theorem \ref{th:R0} is false. That is, if any root from ${\cal R}_0$ of the polynomial $P(x)$ is also a root of the polynomial $h(x)\in R_0[x]$,
dividing $P(x)$ on the right, then it is not necessary that $h(x)$ is a polynomial of the greatest degree from ${\cal R}_0[x]$ dividing $P(x)$ on the right. For example, let $P(x)=(x^2+1)x\in \mathbb {H}[x]$, where
$\mathbb {H}$ is the algebra of Hamiltonian quaternions.
Then any root of the polynomial $P(x)$ from $\mathbb{R}$ is a root of the polynomial $h(x)=x$. But $h(x)$ is not the polynomial of greatest degree in $\mathbb{R}[x]$ that right-divides $P(x)$.
\end{remark}
%%%%%%%%%%%%%%%%%%%%%%%%%%%%%%%%%%%%%%%%%%%%%%%%%%%%%%%%%%%%%%%%%%%%%%%%%%%%%%%%%%%%%%%%%%%%%%%%%%%%%%%%%%%%%%

As a corollary in the notation of the Theorem \ref{th:R0} we obtain

\begin{theorem} \label{th:factorization}
Every polynomial $P(x)\in {\cal R}[x]$ can be uniquely represented as
$$
P(x)=c G(x) H(x),
$$
where $c\in R^*$ is the leading coefficient of the polynomial $P(x)$, $H(x)$ is a monic polynomial with coefficients in the division subring ${\cal R}_0$,
$G(x) \in {\cal R}[x]$ is a monic polynomial that does not have right non-constant divisors from ${\cal R}_0[x]$. Moreover, $H(x)$ is the greatest
common right divisor in ${\cal R}_0[x]$ of the polynomials $b_i(x)$, $i\in I$.
\end{theorem}

%%%%%%%%%%%%%%%%%%%%%%%%%%%%%%%%%%%%%%%%%%%%%%%%%%%%%%%%%%%%%%%%%%%%%%%%%%%%%%%%%%%%%%%%%%%%%%%%%%%%%%%%%%%%%%%

\begin{theorem} \label{th:roots}
Let ${\cal R}$ be a division ring, $P(x)= \sum_{i=0}^n a_ix^i\in {\cal R}[x]$. Suppose that $c$ is not a central element and $P(c)\ne 0$.
Then the conjugacy class $[c]$ has infinitely many elements that are not roots of $P(x)$.
\end{theorem}

%%%%%%%%%%%%%%%%%%%%%%%%%%%%%%%%%%%%%%%%%%%%%%%%%%%%%%%%%%%%%%%%%%%%%%%%%%%%%%%%%%%%%%%%%%%%%%%%%%%%%%%%

Proof. Herstein's theorem (\cite[Theorem 13.26]{La91}) says that the set $[c]$ is infinite. Thus, if $P(x)$ either has no roots
in $[c]$ or has finitely many roots in $[c]$, then the conjugacy class of $[c]$ has infinitely many elements that are not roots of $P(x)$.
Suppose that $P(x)$ has infinitely many roots in $[c]$. It follows from \cite[Proposition 2]{BrWh83} that the set of all $y\in {\cal R}^*$ such that $P(ycy^{-1})=0$
coincides with the set
$$
V := \{ y\in {\cal R}^* | \sum_{i=0}^n a_iyc^i =0\}.
$$
Then $V$ is an infinite set, since $P(x)$ has infinitely many roots in the conjugacy class $[c]$.
Note that $V \cup \{0\}$ is a right vector space over the division ring ${\cal R}^{(c)}$. Since $c$ is not a root of $P(x)$,
we have $1 \notin V$. Then $1+y \notin V$ for any $y \in V$, hence $(1+y)c(1+y)^{-1}$ is not a root of $P(x)$ for any $y \in V$.

Let $y_1, y_2 \in V$, $y_1\ne y_2$. Let us show that
$$
(1+y_1)c(1+y_1)^{-1} \ne (1+y_2)c(1+y_2)^{-1}.
$$ Indeed, if
$(1+y_1)c(1+y_1)^{-1} = (1+y_2)c(1+y_2)^{-1}$, then $1+y_1=(1+y_2)z$ for some $z \in R^{(c)}$. Whence
$$
1=(y_2z-y_1)(1-z)^{-1}.
$$
Which contradicts the fact that $1 \notin V$.

Thus, we obtain infinitely many different elements of the form $(1+y)c(1+y)^{-1}$ that belong to $[c]$ and are not roots of the polynomial $P(x)$. \qed

%%%%%%%%%%%%%%%%%%%%%%%%%%%%%%%%%%%%%%%%%%%%%%%%%%%%%%%%%%%%%%%%%%%%%%%%%%%%%%%%%%%%%%%%%%%%%%%%%%%%%%%%%%%%%%%%%%%%%%%%%%%%%%%%%%%%%%%
%%%%%%%%%%%%%%%%%%%%%%%%%%%%%%%%%%%%%%%%%%%%%%%%%%%%%%%%%%%%%%%%%%%%%%%%%%%%%%%%%%%%%%%%%%%%%%%%%%%%%%%%%%%%%%%%%%%%%%%%%%%%%%%%%%%%%%%

\section{Spherical roots of polynomials over generalized quaternion division algebras}

Let $\cal Q$ be a generalized quaternion division algebra over a field $F$. We will need the following
special case of lemma \cite[Lemma 16.17]{La91}.

%\begin{lemma} [Lem!!!].
%Let $A$ be a division ring, and let $P(x)=L(x)R(x)\in A[x]$. And let $d \in A$ be such that $h:=R(d)\ne 0$.
%Then
%$$
%P(d)=L(hdh^{-1})R(d).
%$$
%In particular, if $d$ is a root of $P(x)$ but not a root of ${\cal R}(x)$, then $hdh^{-1}$ is a root of $L(x)$.
%\end{lemma}

%Thus, the root of the polynomial $d$ is either the root of its right factor, or the element conjugate to $d$ is the root of its left factor.

%%%%%%%%%%%%%%%%%%%%%%%%%%%%%%%%%%%%%%%%%%%%%%%%%%%%%%%%%%%%%%%%%%%%%%%%%%%%%%%%%%%%%%%%%%%%%%%%%%%%%%%%%%%%%%%%%%%%%%%%%%%%%%%%%%%%%%%%%%%%%%%%%%%%%%%%

\begin{lemma} \label{lm:Lam}
Let $\cal Q$ be a generalized quaternion division algebra with center $F$, and let $B$ be a conjugacy class of $\cal Q$ with minimal polynomial $\lambda(x)$ over $F$.
If $P(x)\in {\cal Q}[x]$ has two roots in $B$, then $P(x)\in {\cal Q}[x]\lambda(x)$ and $P(x)$ vanishes on any element of $B$.
\end{lemma}

%%%%%%%%%%%%%%%%%%%%%%%%%%%%%%%%%%%%%%%%%%%%%%%%%%%%%%%%%%%%%%%%%%%%%%%%%%%%%%%%%%%%%%%%%%%%%%%%%%%%%%%%%%%%%%%%%%%%%%%%%%%%%%%%%%%%%%%%%%%%%%%%%%%%%%%%

Recall that a root $q$ of a polynomial $P(x)$ is called spherical if $q$ does not belong to the center of the algebra and any element $d \in [q]$ is also a root of the polynomial $P(x)$.
As a consequence of the lemma \ref{lm:Lam} we obtain

%%%%%%%%%%%%%%%%%%%%%%%%%%%%%%%%%%%%%%%%%%%%%%%%%%%%%%%%%%%%%%%%%%%%%%%%%%%%%%%%%%%%%%%%%%%%%%%%%%%%%%%%%%%%%%%%%%%%%%%%%%%%%%%%%%%%%%%%%%%%%%%%%%%%%%%%

\begin{lemma} \label{le:f1f2}
If a polynomial $P(x)$ of the form (\ref{eq_P(x)}) with coefficients in $\cal Q$ has spherical roots $a_1,\dots,a_n$ lying in different conjugacy classes, then $P(x)$ is divisible by the product of minimal polynomials of these roots.
\end{lemma}

Proof. We prove by induction on the number of spherical root classes. Let $f_i(x)$ be the minimal polynomial of $a_i$, $1\le i \le n$.
It follows from Lemma \ref{lm:Lam} that $P(x)$ is divisible by $f_1(x)$. Suppose that the statement is true for $k$ roots. Then
$$
P(x)=P_1(x)f_k(x)\dots f_1(x).
$$
for some $P_1(x) \in {\cal Q}[x]$.

We prove the statement for $k+1$ roots.
Since $f_k(x)\dots f_1(x) \in F[x]$, then
$$
P(b)=P_1(b)f_k(b)\dots f_1(b)
$$
for any $b\in {\cal Q}$.
Any element from the class $[a_{k+1}]$ is a root of the polynomial $P(x)$, but is not a root of the polynomial $f_k(x)\dots f_1(x)$, since $a_1,\dots,a_{k+1}$ lie in different conjugacy classes.
Then $a_{k+1}$ is a spherical root of the polynomial $P_1(x)$. It follows from Lemma \ref{lm:Lam} that $P_1(x)=P_2(x)f_{k+1}(x)$ for some polynomial $P_2(x)\in {\cal Q}[x]$. Whence
$$
P(x)=P_2(x)f_{k+1}(x)f_k(x)\dots f_1(x).
$$
\qed

%%%%%%%%%%%%%%%%%%%%%%%%%%%%%%%%%%%%%%%%%%%%%%%%%%%%%%%%%%%%%%%%%%%%%%%%%%%%%%%%%%%%%%%%%%%%%%%%%%%%%%%%%%%%%%%%%%%%%%%%%%%%%%%%%%%%%%%%%%%%%%%%%%%%%%%%%%%%%%%%%%%%%

Next, we obtain an estimate for the number of different conjugacy classes of spherical roots depending on the degree of the polynomial.

\begin{theorem} \label{th:sfer_korni}
The polynomial
$$
P(x)=x^n+a_{n-1}x^{n-1}+\dots + a_1x+a_0
$$
with coefficients in the generalized quaternion division algebra $\cal Q$ over the center $F$ has at most $n/2$ different conjugacy classes of spherical roots.
If the degree of the polynomial $P(x)$ is even and there are $n/2$ different  conjugacy  classes of spherical roots, then all the coefficients of $P(x)$ lie in $F$.
If the degree $n$ is odd and there are $(n-1)/2$ different  conjugacy  classes of spherical roots, then all the coefficients  of $P(x)$ lie in the same subfield of the algebra $\cal Q$.
\end{theorem}

Proof. By Lemma \ref{le:f1f2}, the polynomial $P(x)$ is divisible by the product of the minimal polynomials of its spherical roots.
Since the degree of each such minimal polynomial is 2, the maximum number of factors in such a product is $n/2$. Therefore, the number of different conjugacy classes of spherical roots cannot be greater than $n/2$.

Let the degree of $P(x)$ be even. If the number of different conjugacy classes of spherical roots is exactly $n/2$, then $P(x)$ is equal to the product of polynomials from $F[x]$,
and hence the coefficients of $P(x)$ lie in $F$.

Now let the degree of $P(x)$ be odd. If we assume that there are exactly $(n-1)/2$ distinct classes of spherical roots, we get that $P(x)$ has the form $P(x)=(x-a)f(x)$,
where $f(x)\in F[x]$ is the product of minimal polynomials of the spherical roots
of the polynomial $P(x)$, $a\in {\cal Q}$. Then the coefficients of the polynomial $P(x)$ lie in the subfield $F(a)$.\qed

%%%%%%%%%%%%%%%%%%%%%%%%%%%%%%%%%%%%%%%%%%%%%%%%%%%%%%%%%%%%%%%%%%%%%%%%%%%%%%%%%%%%%%%%%%%%%%%%%%%%%%%%%%%%%%%%%%%%%%%%%%%%%%%%%%%%%%%%%%%%%%%%%%%%%%%%

For polynomials of the third degree we obtain a simple sufficient condition for the absence of spherical roots.

%%%%%%%%%%%%%%%%%%%%%%%%%%%%%%%%%%%%%%%%%%%%%%%%%%%%%%%%%%%%%%%%%%%%%%%%%%%%%%%%%%%%%%%%%%%%%%%%%%%%%%%%%%%%%%%%%%%%%%%%%%%%%%%%%%%%%%%%%%%%%%%%%%%%%%%%

\begin{corollary} \label{cor_sfer}
If the coefficients of a polynomial $P(x)=x^3+ax^2+bx+c \in {\cal Q}[x]$ do not lie in the same subfield of $\cal Q$, then the polynomial $P(x)$ has no spherical roots.
\end{corollary}

%Proof. Let $P(x)$ have a spherical root. Then, by Lemma \ref{lm:Lam}, it has the form $P(x)=(x-a)\lambda(x)$, where $\lambda(x)$ is a second-degree polynomial over $F$ (the minimal polynomial of this spherical root). It follows %that $P(x)\in F(a)[x]$, i.e.
%all coefficients of $P(x)$ lie in the subfield $F(a)$. \qed

%%%%%%%%%%%%%%%%%%%%%%%%%%%%%%%%%%%%%%%%%%%%%%%%%%%%%%%%%%%%%%%%%%%%%%%%%%%%%%%%%%%%%%%%%%%%%%%%%%%%%%%%%%%%%%%%%%%%%%%%%%%%%%%%%%%%%%%%%%%%%%%%%%%%%%%%
%%%%%%%%%%%%%%%%%%%%%%%%%%%%%%%%%%%%%%%%%%%%%%%%%%%%%%%%%%%%%%%%%%%%%%%%%%%%%%%%%%%%%%%%%%%%%%%%%%%%%%%%%%%%%%%%%%%%%%%%%%%%%%%%%%%%%%%%%%%%%%%%%%%%%%%%

\section{Spherical roots of polynomials over the algebra of Hamiltonian quaternions}

It is known (see \cite[Theorem 16.14]{La91}) that every non-constant polynomial in one variable with coefficients in the algebra of Hamiltonian quaternions $\mathbb{H}$ has a root in $\mathbb{H}$. Moreover,
a root is either isolated or spherical.
In this section, we apply the results of the previous sections to analyze the existence of spherical roots of some types of polynomials with coefficients in $\mathbb{H}$.

\begin{lemma} \label{lm:sfer}
If $x_1$ is a spherical root of a polynomial $P(x)\in \mathbb{H}[x]$, then each maximal subfield of $\mathbb{H}$ contains a root of $P(x)$ from the conjugacy class $[x_1]$.
\end{lemma}

Proof. The element $x_1$ lies in the maximal subfield $\mathbb{R}(x_1)$ of $\mathbb{H}$.
All maximal subfields of $\mathbb{H}$ are isomorphic to the field of complex numbers $\mathbb{C}$. Then by the Skolem-Noether theorem (\cite[\S 12.6]{Pi86}) all maximal subfields are conjugate, i.e.,
if $K$ is a maximal subfield, then $K=g \mathbb{R}(x_1) g^{-1}$ for some $g\in{\mathbb{H}}, g\ne{0}$. Then $gx_1g^{-1}$ is a root of the polynomial $P(x)$ from the class $[x_1]$ that lies in $K$. \qed

%%%%%%%%%%%%%%%%%%%%%%%%%%%%%%%%%%%%%%%%%%%%%%%%%%%%%%%%%%%%%%%%%%%%%%%%%%%%%%%%%%%%%%%%%%%%%%%%%%%%%%%%%%%%%%%%%%%

\begin{remark}
If a polynomial $P(x)\in \mathbb{H}[x]$ has no roots in some maximal subfield $K\subset \mathbb{H}$, then by Lemma \ref{lm:sfer} this polynomial has no spherical roots.
A similar approach to analyzing the existence of spherical roots of polynomials can be used in the case of quaternion division algebras, in which there are only a finite number of isomorphism classes of maximal subfields. For example, in the case of quaternion algebras over local fields.
\end{remark}

%%%%%%%%%%%%%%%%%%%%%%%%%%%%%%%%%%%%%%%%%%%%%%%%%%%%%%%%%%%%%%%%%%%%%%%%%%%%%%%%%%%%%%%%%%%%%%%%%%%%%%%%%%%%%%%%%%%%

In general, finding the roots of polynomials with coefficients in the algebra of Hamiltonian quaternions is a difficult problem
(see, for example, \cite{FaMiSeSo17}, \cite{HuSo02}, \cite{JaOp10}, \cite{SePeVi01}, \cite{SeSi01}), but for polynomials of a special form one can easily obtain an estimate for the number of  conjugacy  classes of spherical roots.

%%%%%%%%%%%%%%%%%%%%%%%%%%%%%%%%%%%%%%%%%%%%%%%%%%%%%%%%%%%%%%%%%%%%%%%%%%%%%%%%%%%%%%%%%%%%%%%%%%%%%%%%%%%%%%%%%%%%%%%%%%%

\begin{theorem} \label{th:H1}
Consider a polynomial
$$
P(x)=x^n+a_{n-1}x^{n-1}+...+a_kx^k+...+a_mx^m+...+a_1x+a_0 \in \mathbb{H}[x],
$$
such that all its coefficients, except possibly two $a_k$ and $a_m$, lie in $\mathbb{R}$ (the case $m=0$ is possible). Then

1. if one of the coefficients $a_k$ or $a_m$ is from $\mathbb{R}$, and the other is not, then the polynomial $P(x)$ has no spherical roots;

2. if both coefficients $a_k$ and $a_m$ are not from $\mathbb{R}$, but lie in the same subfield, then the polynomial $P(x)$ has at most $(k-m)/2$
different conjugacy classes of spherical roots (in particular, if $k=m+1$, then $P(x)$ has no spherical roots);

3. if both coefficients $a_k$ and $a_m$ do not lie in the same subfield, then the polynomial $P(x)$ has no spherical roots.
\end{theorem}

Proof. 1. Let $a_m \notin{\mathbb{R}}$. Find the roots of the polynomial $P(x)$ that lie in a maximal subfield $K$, $K \ne \mathbb{R}(a_m)$.
As a basis of the algebra $\mathbb{H}$ over the field $K$ we take 1 and $a_m$. Then
$$
P(x)=1\cdot(x^n+a_{n-1}x^{n-1}+...+a_kx^k+...+a_{m-1}x^{m-1}+a_{m+1}x^{m+1}+...+a_1x+a_0)+a_m\cdot x^m
$$
$$
=1\cdot b_1(x)+a_m\cdot b_2(x),
$$
where $b_1(x)$ and $b_2(x)$ are from $\mathbb{R}[x]\subset K[x]$. Since $b_2(x)=x^m$ either has no roots (if $m=0$)
or has a unique root equal to 0 (if $m\ne 0$), then the polynomial $P(x)$, according to Theorem \ref{th:R0},
has no other roots in the field $K$, except possibly 0. Consequently, by Lemma \ref{lm:sfer} in this case $P(x)$ has no spherical roots.
Similar reasoning if $a_k \notin{\mathbb{R}}$.

2. Find the roots of the polynomial $P(x)$ in a maximal subfield $K \ne \mathbb{R}(a_m)$. As a basis of $\mathbb{H}$ over $K$
we take 1 and $a_m$. Since $a_k \in \mathbb{R}(a_m)$, then $a_k=u+va_m$, where $u, v \in \mathbb{R}$. Then, when expanding the polynomial $P(x)$ in basis 1 and $a_m$, we obtain
$$
b_1(x)=x^n+a_{n-1}x^{n-1}+...+a_{k+1}x^{k+1}+ux^k+a_{k-1}x^{k-1}+...+a_{m+1}x^{m+1}+a_{m-1}x^{m-1}+...+a_1x+a_0,
$$
$b_2(x)=x^m+vx^k$. We have, $b_2(x)=x^m+vx^k=x^m(1+vx^{k-m})$. Thus, the roots of the polynomial $P(x)$ in the field $K$ are either 0 or the roots of the polynomial $x^{k-m}+1/v$.
If the polynomial $x^{k-m}+1/v$ has a root in $\mathbb{R}$, then this root is not a spherical root of $P(x)$.
Since $x^{k-m}+1/v \in \mathbb{R}[x]$, then for every root $a$ of this polynomial not in $\mathbb{R}$, the conjugate quaternion $\overline{a}$ is also a root. Since $a$ and $\overline{a}$
belong to the same conjugacy class, then this polynomial has at most $(k-m)/2$ roots that lie in different conjugacy classes.
Thus, $P(x)$ has at most $(k-m)/2$ different conjugacy classes of spherical roots.

3. We will look for roots of the polynomial $P(x)$ in the field $\mathbb{R}(a_m)$. Take basis 1 and $a_k$ of the algebra $\mathbb{H}$ over the field $\mathbb{R}(a_m)$.
Then $P(x)=1\cdot b_1(x)+a_k\cdot b_2(x)$, where $b_2(x)=x^k$ and $b_1(x)$ is some polynomial from $\mathbb{R}(a_m)[x]$.
Since $b_2(x)$ has only the root 0, then in the field $\mathbb{R}(a_m)$ the polynomial $P(x)$ cannot have other roots, except possibly 0.
Then, according to Lemma \ref{lm:sfer} $P(x)$ has no spherical roots in this case. \qed

Next, we consider the case of cubic polynomials.
Note that for quadratic polynomials, explicit formulas for finding roots in the case of the algebra of Hamiltonian quaternions were obtained in \cite{HuSo02}.
%%%%%%%%%%%%%%%%%%%%%%%%%%%%%%%%%%%%%%%%%%%%%%%%%%%%%%%%%%%%%%%%%%%%%%%%%%%%%%%%%%%%%%%%%%%%%%%%%%%%%%%%%%%%%%%%%%%%%%%%%%%

\begin{corollary} \label{cor:sfer_kub}
The polynomial
$$
P(x)=x^3+ax^2+bx+c \in \mathbb{H}[x]
$$
depending on the coefficients has

1. at most one conjugacy class of spherical roots if either a) $a, b, c \in \mathbb{R}$, or b) $a, c \notin \mathbb{R}, b\in \mathbb{R}, a\in\mathbb{R}(c)$, or c) $a, b, c \notin \mathbb{R}, a, c \in \mathbb{R}(b)$;

2. only isolated roots in other cases.
\end{corollary}

%%%%%%%%%%%%%%%%%%%%%%%%%%%%%%%%%%%%%%%%%%%%%%%%%%%%%%%%%%%%%%%%%%%%%%%%%%%%%%%%%%%%%%5

Proof.
1. The estimate for the number of conjugacy classes of spherical roots is obtained from Theorem \ref{th:sfer_korni}.

2. If the conditions on the coefficients of the polynomial from the first point are not satisfied, then the following cases are possible:

a) One coefficient does not lie in the field $\mathbb{R}$, and the remaining coefficients belong to $\mathbb{R}$;

b) $a\in \mathbb{R}$, $b,c \notin \mathbb{R}$, $b \in \mathbb{R}(c)$;

c) $c\in \mathbb{R}$, $a,b \notin \mathbb{R}$, $a \in \mathbb{R}(b)$;

d) the coefficients of the polynomial do not lie in the same subfield of the algebra $\mathbb{H}$.

In case a), the absence of spherical roots follows from point 1 of Theorem \ref{th:H1}.
In cases b) and c), the absence of spherical roots follows from point 2 of Theorem \ref{th:H1}. Finally,
the absence of spherical roots in case d) follows from Corollary \ref{cor_sfer}. \qed

\begin{example}
Let us consider several examples of cubic polynomials with coefficients satisfying the conditions from point 1 of the Corollary \ref{cor:sfer_kub}.

The polynomial $x(x^2-1)$ has no spherical roots, and the polynomial $x(x^2+1)$ has a spherical root $i$.

The polynomial
$$
x^3-ix^2-x+i=(x-i)(x^2-1)
$$
has no spherical roots, and the polynomial
$$
x^3-ix^2+x-i=(x-i)(x^2+1)
$$
has $i$ as a spherical root.

The polynomial
$$
x^3+(2-i)x^2+(1-2i)x-i=(x-i)(x^2+2x+1)=(x-i)(x+1)^2
$$
has no spherical roots, since the polynomial $(x+1)^2$ has roots in $\mathbb{R}$,
and the polynomial
$$
x^3+(1-i)x^2+(1-i)x-i=(x-i)(x^2+x+1)
$$
has a spherical root, this is the root of the polynomial $x^2+x+1$.

Thus, if the conditions from point 1 of the Corollary \ref{cor:sfer_kub} are satisfied, a more subtle reasoning is required to answer the question of the existence of a spherical root for a cubic polynomial.
\end{example}

%%%%%%%%%%%%%%%%%%%%%%%%%%%%%%%%%%%%%%%%%%%%%%%%%%%%%%%%%%%%%%%%%%%%%%%%%%%%%%%%%%%%%%%%%%%%%%%%%%%%%%%

%%%%%%%%%%%%%%%%%%%%%%%%%%%%%%%%%%%%%%%%%%%%%%%%%%%%%%%%%%%%%%%%%%%%%%%%%%%%%%%%%%%%%%%%%%%%%%%%%%%%%%%%%%%%%%

\renewcommand{\refname}


{\normalsize\bf References}
\begin{thebibliography}{00}
\leftskip=-7mm
\parskip=-0mm
\parsep=0mm
\itemsep=0mm
\labelwidth=-12mm

\bibitem{Or33} O. Ore. {\em Theory of non-commutative polynomials}, Ann. of Math. (2),
{\bf 34} (1933), no.~3, 480–-508.
% 1933, vol.\,34, no.~3. pp.\,480–-508.

\bibitem{La91} %%%
T.Y. Lam, {\em A first course in noncommutative rings}, Graduate Texts in Mathematics 131, Springer-Verlag, New York, 1991.

\bibitem{GoMo65} %%%
B. Gordon, T.S. Motzkin. {\em On the zeros of polynomials over division rings}, Trans. Amer.
Math. Soc., {\bf 116} (1965), 218-226.

\bibitem{BrWh83}
U. Bray, G. Whaples. {\em  Polynomials with coefficients from a division ring}, Can. J. Math., {\bf 35} (1983), no.~3,  509--515.

\bibitem{Be79} B. Beck. {\em Sur les \`{e}quations polynomiales dans les quaternions}, Enseign. Math.(2) {\bf 25} (1979), no.~3--4, 193–-201 (1980).

\bibitem{GoTi24} A.G. Goutor, S.V. Tikhonov. {\em Roots of polynomials over division rings}, Dokl. Nats. Akad. Nauk Belarusi, {\bf 68} (2024), no.~5, 359--364.
%https://doi.org/10.29235/1561-8323-2024-68-5-359-364

\bibitem{Pi86}
R.S. Pierce,  {\em Associative algegras}, Stud. Hist. Modern Sci., 9
Grad. Texts in Math., 88, Springer-Verlag, New York-Berlin, 1982.

\bibitem{FaMiSeSo17}
M. I. Falc\~{a}o, F. Miranda, R. Severino, M. J. Soares, {\em Mathematica Tools for Quaternionic Polynomials},
Computational science and its applications.  ICCSA 2017. Part II, 394-408.
Lecture Notes in Comput. Sci., 10405
Springer, Cham, 2017.

\bibitem{HuSo02}
L. Huang,  W. So. {\em Quadratic formulas for quaternions}, Appl. Math. Lett. {\bf 15} (2002), no. 5, 533-540.
%(2003d:12003)

\bibitem{JaOp10} D. Janovsk\'{a}, G. Opfer. {\em A note on the computation of all zeros of simple quaternionic polynomials},
SIAM J. Numer. Anal. {\bf 48} (2010), no. 1, 244-256.

\bibitem{SePeVi01} R. Ser\^{o}dio, E. Pereira, J. Vit\'{o}ria, {\em Computing the zeros of
quaternion polynomials}, Comput. Math. Appl. {\bf 42} (2001), no. 8-9, 1229-1237.
%Numerical methods and computational mechanics
%(Miskolc, 1998). %MR 1851239 (2002f:30061

\bibitem{SeSi01} R. Ser\^{o}dio, L.-S. Siu. {\em Zeros of quaternion polynomials},
Appl. Math. Lett. {\bf 14} (2001),  no. 2, 237-239.
\end{thebibliography}
\end{document}